\documentclass{amsart}
\usepackage{amsfonts}
\usepackage{amsmath}
\usepackage{amssymb}
\usepackage{amsthm}

\newtheorem{theorem}{Theorem}[section]
\newtheorem{prop}{Proposition}[section]
\newtheorem{lemma}{Lemma}[section]
\newtheorem{rem}{Remark}[section]
\newtheorem{cor}{Corollary}[section]
\newtheorem{exmp}{Example}[section]
\newtheorem*{prob}{Problem}

\begin{document}
\author{Mark Pankov}
\title[Characterization of apartments]{Characterization of apartments in polar Grassmannians}
\address{Department of Mathematics and Computer Sciences, University of Warmia and Mazury, Olsztyn, Poland}
\email{pankov@matman.uwm.edu.pl markpankov@gmail.com}

\maketitle
\begin{abstract}
Buildings of types $\textsf{C}_n$ and $\textsf{D}_n$ are defined by rank $n$ polar spaces.
The associated building Grassmannians are polar and half-spin Grassmannians.
Apartments in dual polar spaces and half-spin Grassmannians were characterized in \cite{CKS}.
We characterize apartments in
all polar Grassmannians consisting of non-maximal singular subspaces.
This characterization is a partial case of more general results concerning embeddings of polar Johnson graphs
in polar Grassmann graphs.
\end{abstract}

\section{Introduction}
A {\it building} \cite{Tits} is a simplicial complex $\Delta$ together with a family of subcomplexes
called {\it apartments} and satisfying some axioms.
One of the axioms says that
all apartments are isomorphic to a certain Coxeter complex --- the simplicial complex associated with a Coxeter system.
This Coxeter system defines the type of the building --- if $\textsf{X}$ is the diagram of the Coxeter system then
we say that $\Delta$ is a building of {\it type} $\textsf{X}$.
The vertex set of $\Delta$ can be labeled by the nodes of the diagram $\textsf{X}$.
All vertices corresponding to the same node form a {\it Grassmannian}.
We do not consider more general Grassmannians defined by subsets in the set of nodes \cite{Pasini}.

Let ${\mathcal G}$ be one of the Grassmannians associated with $\Delta$.
The intersection of ${\mathcal G}$ with an apartment is said to be an {\it apartment} in  ${\mathcal G}$.
We say that two distinct vertices $a,b\in {\mathcal G}$ are {\it adjacent} if
there exists a simplex $P\in \Delta$ such that $P\cup\{a\}$ and $P\cup\{b\}$ are chambers
(maximal simplices) of $\Delta$. The Grassmannian ${\mathcal G}$ admits the structure of a partial linear space such that
two distinct points are collinear if and only if they are adjacent elements of ${\mathcal G}$;
this partial linear space is called the {\it Grassmann space} of ${\mathcal G}$.
Let $\Gamma$ be the associated {\it Grassmann graph} (the graph whose vertex set is ${\mathcal G}$ and whose edges are pairs of adjacent vertices).
Let also $\Gamma_{a}$ be the restriction of $\Gamma$ to an apartment of ${\mathcal G}$
(the restrictions of $\Gamma$ to any two apartments of ${\mathcal G}$ are isomorphic).
Consider a few examples:
\begin{enumerate}
\item[(1)] If $\textsf{X}=\textsf{A}_{n-1}$, $n\ge 4$ then $\Delta$ is the flag complex of an $n$-dimensional vector space
and the associated Grassmannians are ${\mathcal G}_{k}(V)$, $k\in \{1,\dots,n-1\}$.
If ${\mathcal G}={\mathcal G}_{k}(V)$ then $\Gamma_{a}$ is isomorphic to the Johnson graph $J(n,k)$.
\item[(2)] If $\textsf{X}=\textsf{C}_{n}$ then  $\Delta$ is the flag complex of a rank $n$ polar space.
If ${\mathcal G}$ is the associated dual polar space
then  $\Gamma_{a}$ is isomorphic to the $n$-dimensional hypercube graph $H_n$.
\item[(3)] If $\textsf{X}=\textsf{D}_{n}$
then  $\Delta$ is the oriflamme complex obtained from a polar space of type $\textsf{D}_{n}$.
If ${\mathcal G}$ is one of the half-spin Grassmannians
then $\Gamma_{a}$ is isomorphic to the $n$-dimensional half-cube graph $\frac{1}{2}H_n$.
\item[(4)] If $\textsf{X}=\textsf{C}_{n}$ or $\textsf{D}_{n}$ and ${\mathcal G}$ is formed by
$k$-dimensional singular subspaces of the associated rank $n$ polar space, $k<n-1$ then
$\Gamma_{a}$ is isomorphic to the polar Johnson graph $PJ(n,k)$ (this graph will be defined in Subsection 2.3).
\end{enumerate}
By \cite{CS},
the image of an embedding of $\Gamma_{a}$ in $\Gamma$ is not necessarily an apartment of ${\mathcal G}$.
In the cases (1) and (3),
apartments of ${\mathcal G}$ can be characterized as the images of embeddings of $\Gamma_{a}$ in $\Gamma$
transferring maximal cliques of $\Gamma_{a}$ to independent subsets of the Grassmann space
\cite[Theorems 2.2 and 4.3]{CKS}.
In the case (2), every maximal clique of $\Gamma_{a}$ consists of two vertices and we use other condition known as
{\it local independence} \cite[Theorem 3.1]{CKS}).
Similar characterizations of apartments in some Grassmannians associated with
buildings of exceptional types were obtained in  \cite{Kasikova3}.
Also, \cite{Kasikova3} contains new proofs of the results from \cite{CKS} mentioned above.
These proofs are based on the following fact:
for any subset ${\mathcal X}\subset{\mathcal G}$ and any $x\in {\mathcal X}$ denote by ${\mathcal X}(x)$
the set formed by $x$ and all elements of ${\mathcal X}$ adjacent with $x$,
then ${\mathcal X}$ is an apartment of ${\mathcal G}$ if
for every $x\in {\mathcal X}$ there is an apartment ${\mathcal A}_{x}\subset {\mathcal G}$, $x\in {\mathcal A}_{x}$
such that ${\mathcal X}(x)={\mathcal A}_{x}(x)$
and some technical conditions hold
(this is a partial case of a more general result \cite[Theorem 3.3]{Kasikova1}).

Some results related with the characterization of apartments as the images of isometric embeddings
can be found in \cite{Pankov1,Pankov2}.

It was noted above that buildings of types $\textsf{C}_n$ and $\textsf{D}_n$ are defined by rank $n$ polar spaces.
The associated building Grassmannians are polar Grassmannians (in particular, dual polar spaces) and half-spin Grassmannians.
Apartments in dual polar spaces and half-spin Grassmannians were characterized in \cite{CKS}.
In this paper we characterize apartments in the remaining polar Grassmannians.

Let $\Pi$ be a polar space of rank $n\ge 3$ and let $\Gamma_{k}(\Pi)$ be the Grassmann graph
formed by $k$-dimensional singular subspaces of $\Pi$.
We investigate embeddings of $PJ(l,m)$ in $\Gamma_{k}(\Pi)$ such that
$$m\le k\le n-2\;\mbox{ and }\;l-m\le n-k.$$
In the case when $l-m=n-k$, we get a characterization of apartments in parabolic subspaces of the associated Grassmann space
(Theorems 4.2 and 4.3). Parabolic subspaces admit the natural structure of polar Grassmann spaces.
As an application, we show that parabolic subspaces can be characterized as convex subspaces isomorphic
to polar Grassmann spaces (Subsection 4.2).

Theorem 4.3 is a partial case of a more general result (Theorem 4.4), but we give an independent proof of Theorem 4.3
based on \cite[Theorem 3.3]{Kasikova1}.
The latter result can not be exploited if $l-m<n-k$. In this case, the images of embeddings of $PJ(l,m)$ in $\Gamma_{k}(\Pi)$
can not be apartments of polar Grassmannians.

\section{Polar Grassmannians}
\subsection{Partial linear space}
Let $P$ be a non-empty set and let ${\mathcal L}$ be a family of proper subsets of $P$.
Elements of $P$ and ${\mathcal L}$ will be called {\it points} and {\it lines}, respectively.
We say that two or more points are {\it collinear} if there is a line containing all of them.
The pair $\Pi=(P,{\mathcal L})$ is a {\it partial linear space} if
the following axioms hold:
\begin{enumerate}
\item[$\bullet$]
every line contains at least two points and every point belongs to a line;
\item[$\bullet$]
for any distinct collinear points $p,q\in P$ there is precisely one line containing them,
this line will be denoted by $p\,q$.
\end{enumerate}
Let $\Pi=(P,{\mathcal L})$ be a partial linear space.
We say that $S\subset P$ is a {\it subspace} of $\Pi$
if for any distinct collinear points $p,q\in S$
the line $p\,q$ is contained in $S$.
A subspace is called {\it singular} if any two distinct points of the subspace are collinear.
The empty set and a single point are singular subspaces.
Using Zorn lemma, we show that every singular subspace is contained in
a maximal singular subspace.
For every subset $X\subset P$ the minimal subspace containing $X$,
i.e. the intersection of all subspaces containing $X$, is called {\it spanned} by $X$
and denoted by $\langle X\rangle$.
We say that $X$ is an {\it independent} subset if the subspace $\langle X\rangle$ can not be spanned by a proper subset of $X$.

Let $S$ be a subspace of $\Pi$ (possible $S = P$).
An independent subset $X\subset S$ is said to be a {\it base} of $S$ if $\langle X\rangle=S$.
The {\it dimension} of $S$ is defined as the smallest cardinality $\alpha$
such that $S$ has a base of cardinality $\alpha+1$.
The dimension of the empty set and a single point is equal to $-1$ and $0$ (respectively),
lines are $1$-dimensional subspaces,
$2$-dimensional singular subspaces are called {\it planes}.

The {\it collinearity graph} $\Gamma_{\Pi}$ is the graph whose vertex set is $P$ and whose edges are
pairs of distinct collinear points.
Suppose that $\Gamma_{\Pi}$ is connected and define the {\it distance} $d(p,q)$ between points $p,q\in P$
as the smallest number $i$ such that there is a path of the length $i$ connecting $p$ and $q$;
a path connecting $p$ and $q$ is said to be a {\it geodesic} if it consists of $d(p,q)$ edges.
A subspace $S$ is called {\it convex} if for any points $p,q\in S$
every geodesic connecting $p$ and $q$ is contained in $S$.
For any $X\subset P$ the intersection of all convex subspaces containing $X$
is said to be the {\it convex closure} of $X$.

Two partial linear spaces $\Pi=(P,{\mathcal L})$ and $\Pi'=(P',{\mathcal L}')$ are {\it isomorphic}
if there is a bijection $f:P\to P'$ such that $f({\mathcal L})={\mathcal L}'$.
Any such bijection is called a {\it collineation} of $\Pi$ to $\Pi'$.
Every collineation transfers subspaces to subspaces and induces an isomorphism between the collinearity graphs;
moreover, it is distance preserving and the images of convex subspaces are convex subspaces
(if the distance is well-defined).

\subsection{Polar spaces and polar Grassmannians}
A {\it polar space} is a partial linear space $\Pi=(P,{\mathcal L})$
satisfying the following axioms:
\begin{enumerate}
\item[$\bullet$] every line contains at least three points,
\item[$\bullet$] there is no point collinear with all other points,
\item[$\bullet$] if $p\in P$ and $L\in {\mathcal L}$
then $p$ is collinear with precisely one or all points of $L$,
\item[$\bullet$] any flag formed by singular subspaces is finite.
\end{enumerate}
If a polar space $\Pi=(P,{\mathcal L})$ contains a singular subspace whose dimension is greater than $1$
then all maximal singular subspaces of $\Pi$ are projective spaces
of the same dimension $m\ge 2$ and the number $m+1$ is called the {\it rank} of the polar space.
The collinearity relation of $\Pi$ is denoted by $\perp$:
for points $p,q\in P$ we write $p\perp q$ if $p$ is collinear with $q$
and $p\not\perp q$ otherwise.
If $X,Y\subset P$ then $X\perp Y$ means that every point of $X$
is collinear with all points of $Y$.
For any $X\subset P$ we denote by $X^{\perp}$ the set of all points
$p\in P$ such that $p\perp X$.

Let $\Pi=(P,{\mathcal L})$ be a polar space of rank $n$.
Denote by ${\mathcal G}_k(\Pi)$ the Grassmannian formed by $k$-dimensional singular subspaces of $\Pi$.
Then ${\mathcal G}_0(\Pi)=P$ and ${\mathcal G}_{n-1}(\Pi)$ consists of maximal singular subspaces.
Two elements of ${\mathcal G}_{n-1}(\Pi)$ are {\it adjacent} if their intersection is $(n-2)$-dimensional.
In the case when $k\le n-2$, elements $S,U\in{\mathcal G}_k(\Pi)$ are {\it adjacent} if
$$\dim(S\cap U)=k-1\;\mbox{ and }\;S\perp U.$$
The associated {\it Grassmann graph} is denoted by $\Gamma_{k}(\Pi)$.
This graph is connected.

Let $S$ and $U$ be incident singular subspaces of $\Pi$ such that
$\dim S<k<\dim U$.
Denote by $[S,U]_k$ the set of all $X\in {\mathcal G}_k(\Pi)$ satisfying
$S\subset X\subset U$.
In the case when $S=\emptyset$, this set will be denoted by $\langle U]_k$.
Also, we  write $[S\rangle_k$ for the set formed by all elements of ${\mathcal G}_k(\Pi)$
containing $S$.

The {\it dual polar space} ${\mathfrak G}_{n-1}(\Pi)$ is the partial linear space
whose points are elements of ${\mathcal G}_{n-1}(\Pi)$ and whose lines are subsets of type
$$[S\rangle_{n-1},\;S\in {\mathcal G}_{n-2}(\Pi).$$
The {\it polar Grassmann space} ${\mathfrak G}_k(\Pi)$, $0\le k\le n-2$ is the partial linear space
whose points are elements of ${\mathcal G}_k(\Pi)$ and whose lines are subsets of type
$$[S,U]_k,\;\;S\in {\mathcal G}_{k-1}(\Pi),\;U\in {\mathcal G}_{k+1}(\Pi).$$
Note that ${\mathfrak G}_0(\Pi)=\Pi$ and $\Gamma_{k}(\Pi)$ is the collinearity graph of ${\mathfrak G}_k(\Pi)$.

\begin{lemma}\cite[Proposition 4.15]{Pankov-book}\label{lemma2-1}
If an element of ${\mathcal G}_{k}(\Pi)$ is collinear with two distinct points on a line of ${\mathfrak G}_k(\Pi)$
then it is collinear with all points on this line.
\end{lemma}

The class of maximal singular subspaces of ${\mathfrak G}_k(\Pi)$ coincides with the class of maximal cliques of $\Gamma_{k}(\Pi)$.
Every maximal singular subspace of ${\mathfrak G}_{n-1}(\Pi)$ is a line.
In the case when $1\le k\le n-2$, there are precisely the following two types of maximal singular subspaces of ${\mathfrak G}_k(\Pi)$
\cite[Proposition 4.16]{Pankov-book}:
\begin{enumerate}
\item[$\bullet$] the {\it top} $\langle U]_{k}$, $U\in {\mathcal G}_{k+1}(\Pi)$,
\item[$\bullet$] the {\it star} $[S,M]_{k}$,  $S\in {\mathcal G}_{k-1}(\Pi)$, $M\in {\mathcal G}_{n-1}(\Pi)$ and $k\le n-3$.
\end{enumerate}
Tops and stars are projective spaces of dimension $k+1$ and $n-k-1$, respectively.

Every {\it big star}  $[S\rangle_{k}$, $S\in {\mathcal G}_{k-1}(\Pi)$ is a polar space of rank $n-k$ \cite[Lemma 4.4]{Pankov-book}.

Let $N$ be a $(k-m-1)$-dimensional singular subspace of $\Pi$.
Then $[N\rangle_{k}$ is a subspace of ${\mathfrak G}_{k}(\Pi)$;
subspaces of such type are called {\it parabolic}.
If $S\in [N\rangle_{k}$ then $[N,S]_{k-m}$ is an $m$-dimensional singular subspace of
the polar space $[N\rangle_{k-m}$. This correspondence is a collineation
of the parabolic subspace $[N\rangle_{k}$ to the index $m$ Grassmann space of $[N\rangle_{k-m}$.

\subsection{Apartments in polar Grassmannians}
A {\it frame} of $\Pi$ is a set consisting of $2n$ distinct points $p_1,\dots,p_{2n}\in P$
such that for each $i\in \{1,\dots,2n\}$ there is unique $\sigma(i)\in \{1,\dots,2n\}$
satisfying $p_i\not\perp p_{\sigma(i)}$.
Every frame is an independent subset of $\Pi$.
Thus any $k+1$ mutually collinear points of a frame span a $k$-dimensional singular subspace.

Let $B=\{p_{1},\dots,p_{2n}\}$ be a frame of $\Pi$.
The associated apartment ${\mathcal A}_{k}\subset {\mathcal G}_{k}(\Pi)$
is formed by all $k$-dimensional singular subspaces spanned by subsets of $B$,
i.e. it consists of all subspaces
$$\langle p_{i_1},\dots,p_{i_{k+1}}\rangle$$
such that
$$\{i_{1},\dots,i_{k+1}\}\cap \{\sigma(i_{1}),\dots,\sigma(i_{k+1})\}=\emptyset.$$
Note that ${\mathcal A}_{0}=B$.

For every subset ${\mathcal X}\subset {\mathcal G}_k(\Pi)$ we denote by $\Gamma({\mathcal X})$ the restriction of $\Gamma_{k}(\Pi)$ to ${\mathcal X}$.
Then $\Gamma({\mathcal A}_{n-1})$ is isomorphic to the $n$-dimensional hypercube graph $H_n$,
every maximal clique of this graph is a pair of adjacent vertices.
If $1\le k\le n-2$ then  there are precisely the following two types of maximal cliques of $\Gamma({\mathcal A}_{k})$:
\begin{enumerate}
\item[$\bullet$] the {\it top} ${\mathcal T}_{k}(U):={\mathcal A}_{k}\cap\langle U]_{k}$, $U\in {\mathcal A}_{k+1}$,
\item[$\bullet$] the {\it star} ${\mathcal S}_{k}(S,M):={\mathcal A}_{k}\cap[S,M]_{k}$,  $S\in {\mathcal A}_{k-1}$, $M\in {\mathcal A}_{n-1}$ and $k\le n-3$.
\end{enumerate}
Each of these cliques is a base of the maximal singular subspace containing it.
For every $S\in {\mathcal A}_{k-1}$ the {\it big star} ${\mathcal B}_{k}(S)$ is the intersection of
${\mathcal A}_{k}$ with the big star $[S\rangle_{k}$; this is a frame of the polar space $[S\rangle_{k}$.

Consider the set
$$J:=\{1,\dots,n,-1,\dots,-n\}.$$
We say that a subset $X\subset J$ is {\it singular} if
$$j\in X\;\Longrightarrow\; -j\not\in X.$$
We define the {\it polar Johnson graph} $PJ(n,k)$, $k\in\{0,1,\dots,n-1\}$ as the graph whose vertex set is formed by all
singular subsets consisting of $k+1$ elements. Two such subsets are adjacent
(connected by an edge) if their intersection consists of $k$ elements,
and in the case when $k\le n-2$, we also require that their sum is
singular. The graphs $PJ(n,k)$ and $\Gamma({\mathcal A}_{k})$ are isomorphic.
A maximal clique of $PJ(n,k)$, $1\le k\le n-2$ is said to be a {\it top} or a {\it star} if
it is defined by a vertex of $PJ(n,k+1)$ or a vertex of $PJ(n,k-1)$ and a vertex of $PJ(n,n-1)$,
respectively.
If $1\le k\le n-3$ then every $X\in PJ(n,k-1)$ defines the {\it big star} of $PJ(n,k)$ which
consists of all vertices of $PJ(n,k)$ containing $X$.
The restriction of $PJ(n,k)$ to every big star is isomorphic to $PJ(n-k,0)$.

\begin{lemma}\label{lemma2-2}
Let $X$ and $Y$ be subsets of $\Pi$ such that
$$|Y|\le|X|\le n,\;X\perp X,\;Y\perp Y$$
and for every point $y\in Y$ there is unique point $x\in X$ such that $x\not\perp y$.
Then $X\cup Y$ can be extended to a frame of $\Pi$.
\end{lemma}

\begin{proof}
Let $x_{1},\dots,x_{u}$ and $y_{1},\dots, y_{v}$, $v\le u\le n$ be the elements of $X$ and $Y$, respectively.
Suppose that $x_{i}\not\perp y_{i}$, for every $i\le v$.
Then
\begin{equation}\label{eq2-1}
x^{\perp}_{1}\cap y^{\perp}_{1}\cap\dots\cap x^{\perp}_{v}\cap y^{\perp}_{v}
\end{equation}
is a polar space of rank $n-v$, see \cite[Lemma 4.3]{Pankov-book}.
This polar space contains $x_{v+1},\dots,x_{u}$ and we choose two disjoint maximal singular subspaces
of \eqref{eq2-1} such that one of them contains $x_{v+1},\dots,x_{u}$.
As in \cite[Subsection 4.1.4]{Pankov-book}, we show that $x_{v+1},\dots,x_{u}$
can be extended to a frame $B'$ of \eqref{eq2-1}.
Then $B'\cup\{x_{1},y_{1},\dots,x_{v},y_{v}\}$
is a frame of $\Pi$ containing $X\cup Y$.
\end{proof}

Let $l$ be a natural number not greater than $n$.
An $l$-{\it frame} of $\Pi$ is a set consisting of $2l$ distinct points $p_1,\dots,p_{2l}\in P$
such that for each $i\in \{1,\dots,2l\}$ there is unique $\sigma(i)\in \{1,\dots,2l\}$
satisfying $p_i\not\perp p_{\sigma(i)}$. This is a frame of $\Pi$ if $l=n$.
In the general case, every $l$-frame can be extended to a frame (Lemma \ref{lemma2-2}).

Let $B$ be an $l$-frame of $\Pi$ and let $k\le l-1$.
Denote by ${\mathcal X}$  the set formed by all elements of ${\mathcal G}_{k}(\Pi)$ spanned by subsets of $B$.
Then ${\mathcal X}$ is the intersection of all apartments of ${\mathcal G}_{k}(\Pi)$ defined by frames containing $B$.
Also, $\Gamma({\mathcal X})$ is isomorphic to $PJ(l,k)$.

\subsection{Apartments in parabolic subspaces}
Let $S\in {\mathcal G}_{k-1}(\Pi)$.
Then $[S\rangle_{k}$ is a polar space.
If ${\mathcal B}$ is a frame of $[S\rangle_{k}$ then there is a frame $B$ of $\Pi$
such that $S$ is spanned by a subset of $B$ and ${\mathcal B}$ is a big star
in the associated apartment of ${\mathcal G}_{k}(\Pi)$ \cite[Lemma 4.4]{Pankov-book}.

Let $N$ be a $(k-m-1)$-dimensional singular subspace of $\Pi$.
Since the parabolic subspace $[N\rangle_{k}$ can be identified with the index $m$ Grassmann space of the polar space
$[N\rangle_{k-m}$,
frames of $[N\rangle_{k-m}$ define {\it apartments} in $[N\rangle_{k}$.
If ${\mathcal A}$ is an apartment of $[N\rangle_{k}$
then there is a frame $B$ of $\Pi$ such that $N$ is spanned by a subset of $B$
and ${\mathcal A}$ is the intersection of $[N\rangle_{k}$ with the associated apartment of ${\mathcal G}_{k}(\Pi)$.
The restriction of $\Gamma_{k}(\Pi)$ to every apartment of $[N\rangle_{k}$ is isomorphic to  $PJ(n-k+m,m)$.

Let ${\mathcal B}$ be an $l$-frame of the polar space $[N\rangle_{k-m}$
(then $l\le n-k+m$). Suppose that $m\le l-1$ and denote by ${\mathcal X}$
the set formed by all elements of $[N\rangle_{k}$ spanned by elements of ${\mathcal B}$.
Then $\Gamma({\mathcal X})$ is isomorphic to $PJ(l,m)$.

\section{Embeddings of $PJ(l,m)$ in $\Gamma_{k}(\Pi)$}
\setcounter{equation}{0}
Recall that an {\it embedding} of a graph $\Gamma$ in a graph $\Gamma'$
is an injective mapping of the vertex set of $\Gamma$ to the vertex set of $\Gamma'$
such that two vertices of $\Gamma$ are adjacent if and only if their images are adjacent.

In this section we establish some simple facts concerning embeddings of polar Johnson graphs in
the Grassmann graph $\Gamma_{k}(\Pi)$.
First of all, we describe embeddings of  $PJ(l,0)$ in $\Gamma_{k}(\Pi)$ for $l\ge 3$.

\begin{exmp}{\rm
Every big star $[S\rangle_{k}$, $S\in{\mathcal G}_{k-1}(\Pi)$ is a polar space of rank $n-k$.
If ${\mathcal X}$ is an $l$-frame of $[S\rangle_{k}$, $l\le n-k$ then $\Gamma({\mathcal X})$ is isomorphic to $P(l,0)$.
}\end{exmp}

\begin{exmp}{\rm
Suppose that $1\le k\le n-3$ and $N\in {\mathcal G}_{k-2}(\Pi)$, $M\in {\mathcal G}_{k+2}(\Pi)$ are incident.
Then $[N,M]_{k}$ is a polar space of rank $3$
(it is well-known that the Grassmann space formed by lines of a $3$-dimensional projective space is a polar space of rank $3$).
The restriction of $\Gamma_{k}(\Pi)$ to every frame of $[N,M]_{k}$ is isomorphic to $PJ(3,0)$.
}\end{exmp}

\begin{prop}\label{prop3-1}
Let ${\mathcal X}$ be a subset of ${\mathcal G}_{k}(\Pi)$ such that $\Gamma({\mathcal X})$ is isomorphic to $PJ(l,0)$ and $l\ge 3$.
Then $l\le n-k$ and one of the following possibilities
is realized:
\begin{enumerate}
\item[{\rm (1)}] there exists $S\in {\mathcal G}_{k-1}(\Pi)$ such that ${\mathcal X}$ is an $l$-frame of
$[S\rangle_{k}$ {\rm(}if $l=n-k$ then ${\mathcal X}$ is a frame{\rm)},
\item[{\rm (2)}] $l=3$ and there exist incident $N\in {\mathcal G}_{k-2}(\Pi)$ and $M\in {\mathcal G}_{k+2}(\Pi)$
such that ${\mathcal X}$ is a frame of $[N,M]_{k}$.
\end{enumerate}
\end{prop}

\begin{proof}
We take any $X\in {\mathcal X}$. There is unique $Y\in {\mathcal X}$ non-adjacent with
$X$ and all other elements of ${\mathcal X}$ are adjacent with both $X$ and $Y$.
One of the following possibilities is realized:
\begin{enumerate}
\item[{\rm (1)}] $\dim (X\cap Y)=k-1$ and $X\not\perp Y$,
\item[{\rm (2)}] $\dim (X\cap Y)=k-2$.
\end{enumerate}
In the first case, every $Z\in {\mathcal X}\setminus\{X,Y\}$ contains
the $(k-1)$-dimensional singular subspace $S=X\cap Y$
(otherwise $Z$ intersects $X$ and $Y$ in distinct $(k-1)$-dimensional singular subspaces which implies that $X\perp Y$).
So, ${\mathcal X}$ is contained in the big star $[S\rangle_{k}$.
If $l\le n-k$ then ${\mathcal X}$ is an $l$-frame of $[S\rangle_{k}$ (by the definition).

In the case when $l> n-k$, we consider ${\mathcal X}'\subset {\mathcal X}$ such that $\Gamma({\mathcal X}')$ is isomorphic to $PJ(n-k,0)$.
Then ${\mathcal X}'$ is a frame of $[S\rangle_{k}$ and every element of ${\mathcal X}\setminus {\mathcal X}'$
is adjacent with all elements of ${\mathcal X}'$.
Since a polar space does not contain a point collinear with all points of a frame \cite[Corollary 4.1]{Pankov-book},
the inequality $l> n-k$ is impossible.

Consider the case (2). Let $N=X \cap Y$.
An easy verification shows that every element of ${\mathcal X}\setminus\{X,Y\}$
intersects $X$ and $Y$ in $(k-1)$-dimensional singular subspaces containing $N$.
Thus ${\mathcal X}\subset [N\rangle_{k}$ and we can consider elements of ${\mathcal X}$
as lines in the polar space $[N\rangle_{k-1}$.
We take non-adjacent $X',Y'\in{\mathcal X}\setminus\{X,Y\}$.
Then
$$P_{X}=X\cap X',\;P_{Y}=Y\cap X',\;Q_{X}=X\cap Y',\;Q_{Y}=Y\cap Y'$$
are distinct points of $[N\rangle_{k-1}$.
Now consider any non-adjacent
$$X'',Y''\in{\mathcal X}\setminus\{X,Y,X',Y'\}.$$
Each of these lines intersects $X'$ in $P_X$ or $P_{Y}$ and it intersects $Y'$ in $Q_X$ or $Q_{Y}$
(for example, if $X'\cap X''$ does not coincide with $P_X$ or $P_{Y}$ then $X$ and $Y$ both are contained in
the plane spanned by $X'$ and $X''$ which is impossible).
Therefore, one of $X'',Y''$ is $\langle P_{X},Q_{Y}\rangle $ and the other coincides with $\langle P_{Y},Q_{X}\rangle$.
This implies that ${\mathcal X}$ consists of $6$ elements, i.e. $l=3$.
Moreover, $X\perp Y$.
The singular subspace $M=\langle X,Y\rangle $ is $(k+2)$-dimensional and ${\mathcal X}$ is a frame of
$[N,M]_{k}$.
\end{proof}

If $l-m\ge 3$ and $m\ge 1$ then maximal cliques of $PJ(l,m)$ are tops and stars.
Every embedding of $PJ(l,m)$ in $\Gamma_{k}(\Pi)$ transfers them to cliques of $\Gamma_{k}(\Pi)$
--- subsets in stars or tops of ${\mathcal G}_{k}(\Pi)$.
The image of a maximal clique of $PJ(l,m)$ is not necessarily contained in a unique maximal clique of $\Gamma_{k}(\Pi)$.

\begin{prop}\label{prop3-2}
If $l-m\ge 3$ and $m\ge 1$ then the following assertions are fulfilled:
\begin{enumerate}
\item[{\rm (1)}] There are no embeddings of $PJ(l,m)$ in $\Gamma_{k}(\Pi)$ if $l-m>n-k$.
\item[{\rm (2)}] If $l-m\le n-k$ then every embedding of $PJ(l,m)$  in $\Gamma_{k}(\Pi)$ transfers stars of $PJ(l,m)$
to independent subsets of ${\mathfrak G}_{k}(\Pi)$;
moreover, if $l-m\ge 4$ then the image of every star is contained in a star.
\item[{\rm (3)}]
If $4\le l-m\le n-k$ then every embedding of $PJ(l,m)$ in $\Gamma_{k}(\Pi)$ sends every big star of $PJ(l,m)$
to an $(l-m)$-frame in a big star of ${\mathcal G}_{k}(\Pi)$ {\rm(}this is a frame of a big star if $l-m=n-k${\rm)}.
\end{enumerate}
\end{prop}

\begin{proof}
Every big star of $PJ(l,m)$ is isomorphic to $PJ(l-m,0)$.
So, the statements (1) and (3) follow from Proposition \ref{prop3-1}.

Every star of $PJ(l,m)$ is contained in a big star of $PJ(l,m)$.
By Proposition \ref{prop3-1}, its image is contained in a frame of a polar space which is a subspace of ${\mathfrak G}_{k}(\Pi)$.
The image of our star is an independent subset of ${\mathfrak G}_{k}(\Pi)$, since every subset in a frame is independent.
The second part of the statement (2) is a consequence of the statement (3).
\end{proof}

\begin{prop}\label{prop3-3}
If $l\ge 4$ then every embedding of $PJ(l,1)$ in $\Gamma_{k}(\Pi)$ maps all maximal cliques of $PJ(l,1)$
to independent subsets of ${\mathfrak G}_{k}(\Pi)$.
\end{prop}

\begin{proof}
Let $f$ be an embedding of $PJ(l,1)$ in $\Gamma_{k}(\Pi)$ and $l\ge 4$.
By Proposition \ref{prop3-2}, the image of every star is an independent subset of ${\mathfrak G}_{k}(\Pi)$.

Let ${\mathcal T}$ be a top of $PJ(l,1)$ and let ${\mathcal S}$ be a star of $PJ(l,1)$ intersecting ${\mathcal T}$
precisely in two elements.
The top ${\mathcal T}$ consists of three elements which will be denoted by $X_{1},X_{2},X_{3}$.
Suppose that $X_{1},X_{2}\in {\mathcal S}$. If $f({\mathcal T})$ is not independent in ${\mathfrak G}_{k}(\Pi)$
then $f(X_{3})$ is on the line of ${\mathfrak G}_{k}(\Pi)$ joining $f(X_{1})$ and $f(X_{2})$.
By Lemma \ref{lemma2-1}, the latter means that every element of $f(\mathcal S)$ is adjacent with $f(X_{3})$.
This is impossible, since $f$ is an embedding and $X_{3}$ is not adjacent with all elements of ${\mathcal S}$.
\end{proof}

\section{Main results}
\setcounter{equation}{0}
\subsection{Characterization of apartments}
By Proposition \ref{prop3-1}, every subset ${\mathcal X}\subset {\mathcal G}_{k}(\Pi)$
such that $\Gamma({\mathcal X})$ is isomorphic to $PJ(l,0)$, $l\ge 3$
can be extended to a frame of a certain polar space contained in ${\mathfrak G}_{k}(\Pi)$.
On the other hand, there exists a subset ${\mathcal X}\subset {\mathcal G}_{n-1}(\Pi)$ such that
$\Gamma({\mathcal X})$ is isomorphic to $H_n$ and ${\mathcal X}$ is not an apartment of ${\mathcal G}_{n-1}(\Pi)$
\cite[Example 2]{CS}.
So, we need some additional conditions to characterize apartments.

We say that a subset ${\mathcal X}\subset {\mathcal G}_{n-1}(\Pi)$ is {\it locally independent} if for every $S\in {\mathcal X}$
$$\{\;S\cap U\;:\;U\in{\mathcal X} \mbox{ is adjacent with } S\;\}$$
is an independent subset of the projective space $\langle S]_{n-2}$. Apartments of ${\mathcal G}_{n-1}(\Pi)$
and apartments in parabolic subspaces of ${\mathfrak G}_{n-1}(\Pi)$
are locally independent.

\begin{theorem}[B.N. Cooperstein, A. Kasikova, E.E. Shult \cite{CKS}]\label{theorem4-1}
If ${\mathcal X}$ is a locally independent subset of ${\mathcal G}_{n-1}(\Pi)$ and
$\Gamma({\mathcal X})$ is isomorphic to $H_m$ then
there exists an $(n-m-1)$-dimensional singular subspace $N$ such that ${\mathcal X}$ is an apartment in
the parabolic subspace $[N\rangle_{n-1}$.
\end{theorem}

Our first result is the following.

\begin{theorem}\label{theorem4-2}
Let $l$ and $m$ be natural numbers satisfying
\begin{equation}\label{eq4-1}
0<m\le k,\;l\le n\;\mbox{ and }\;l-m=n-k>1.
\end{equation}
Let also ${\mathcal X}$ be a subset of ${\mathcal G}_{k}(\Pi)$ such that
$\Gamma({\mathcal X})$ is isomorphic to $PJ(l,m)$ and
every maximal clique of $\Gamma({\mathcal X})$ is an independent subset of ${\mathfrak G}_{k}(\Pi)$.
Suppose that  one of the following conditions holds:
\begin{enumerate}
\item[${\rm (1)}$]
$2m+2>l$,
\item[${\rm (2)}$]
${\mathcal X}$ is not contained in any big star $[S\rangle_{k}$, $S\in {\mathcal G}_{k-1}(\Pi)$.
\end{enumerate}
Then there exists a $(k-m-1)$-dimensional singular subspace $N$ such that ${\mathcal X}$ is an apartment in
the parabolic subspace $[N\rangle_{k}$.
\end{theorem}

The fact that for a subset ${\mathcal X}\subset {\mathcal G}_{k}(\Pi)$
the graph $\Gamma({\mathcal X})$ is isomorphic to $PJ(l,m)$ and
every maximal clique of $\Gamma({\mathcal X})$ is an independent subset of ${\mathfrak G}_{k}(\Pi)$
can be reformulated in the following form:
${\mathcal X}$ is the image of an embedding of $PJ(l,m)$ in $\Gamma_{k}(\Pi)$
sending maximal cliques of $PJ(n,k)$ to independent subsets of ${\mathfrak G}_{k}(\Pi)$.
If $1\le m\le l-2$ then
maximal cliques of $PJ(l,m)$ are tops and stars
(the second possibility is realized only in the case when $m\le l-3$).
By Propositions \ref{prop3-2} and \ref{prop3-3}, the following assertions are fulfilled:
\begin{enumerate}
\item[$\bullet$] if $l-m\ge 3$ then every embedding of $PJ(l,m)$ in $\Gamma_{k}(\Pi)$ transfers each star of $PJ(l,m)$
to an independent subset of ${\mathfrak G}_{k}(\Pi)$,
\item[$\bullet$] if $l\ge 4$ then every embedding of $PJ(l,1)$ in $\Gamma_{k}(\Pi)$ sends all maximal cliques of $PJ(l,1)$
to independent subsets of ${\mathfrak G}_{k}(\Pi)$.
\end{enumerate}

The second fact implies the following.

\begin{cor}\label{cor4-1}
Suppose that $n-k\ge 3$ and ${\mathcal X}$ is a subset of ${\mathcal G}_{k}(\Pi)$ such that
$\Gamma({\mathcal X})$ is isomorphic to $PJ(n-k+1,1)$.
If ${\mathcal X}$ is not contained in any big star $[S\rangle_{k}$, $S\in {\mathcal G}_{k-1}(\Pi)$
then there exists a $(k-2)$-dimensional singular subspace $N$ such that ${\mathcal X}$ is an apartment in
the parabolic subspace $[N\rangle_{k}$.
\end{cor}

If $m=l-2$ then Theorem \ref{theorem4-2} is a simple consequence of Theorem \ref{theorem4-1}.
The equality $n-k=l-m=2$ implies that $k=n-2$.
Since all maximal cliques of $PJ(l,l-2)$  and $\Gamma_{n-2}(\Pi)$ are tops,
every embedding $f$ of $PJ(l,l-2)$ in $\Gamma_{n-2}(\Pi)$ induces an injective mapping $g$ of the vertex set of
$H_{l}$ to ${\mathcal G}_{n-1}(\Pi)$.
An easy ve\-ri\-fication shows that the image of $g$ is locally independent
if $f$ transfers tops to independent subsets of tops.
By Theorem \ref{theorem4-1}, there is an $(n-l-1)$-dimensional singular subspace $N$
such that the image of $g$ is an apartment in
the parabolic subspace $[N\rangle_{n-1}$. The image of $f$ is the associated apartment in $[N\rangle_{n-2}$.
Since $n-k=l-m$, the dimension of $N$ is equal to $k-m-1$.

To prove Theorem \ref{theorem4-2} in the case when $m\le l-3$, we will use the following result.

\begin{theorem}\label{theorem4-3}
Let $l$ and $m$ be natural numbers satisfying \eqref{eq4-1} and let $m\le l-3$.
If an embedding of $PJ(l,m)$ in $\Gamma_{k}(\Pi)$ sends tops to independent subsets of ${\mathfrak G}_{k}(\Pi)$
contained in tops
then there exists a $(k-m-1)$-dimensional singular subspace $N$ such that ${\mathcal X}$ is an apartment in
the parabolic subspace $[N\rangle_{k}$.
\end{theorem}

\begin{prob}{\rm
Let $l$ and $m$ be natural numbers satisfying \eqref{eq4-1} and let $l\ge 2m+2$,
i.e. the condition (1) from Theorem \ref{theorem4-2} does not hold.
Let also $f$ be an embedding of $PJ(l,m)$ in $\Gamma_{k}(\Pi)$ sending maximal cliques of $PJ(l,m)$
to independent subsets of ${\mathfrak G}_{k}(\Pi)$.
Show that the image of $f$ is not contained in any big star of ${\mathcal G}_{k}(\Pi)$ or construct a contrexample.
Since every big star of ${\mathcal G}_{k}(\Pi)$ is a polar space of rank $n-k=l-m$, our question can be reformulated as follows:
is there an embedding of $PJ(l,m)$ in the collinearity graph of a rank $l-m$ polar space
such that maximal cliques go to independent subsets?
}\end{prob}

\subsection{Application: Convex subspaces of polar Grassmann spaces}
In some cases, the Grassmann space ${\mathfrak G}_{k}(\Pi)$ is not spanned by an apartment of ${\mathcal G}_{k}(\Pi)$
\cite{BB,CS,Pankov-book}.
By \cite{Kasikova2}, the convex closure of every apartment of ${\mathcal G}_{k}(\Pi)$ coincides with ${\mathfrak G}_{k}(\Pi)$;
moreover, if ${\mathcal A}$ is an apartment in a parabolic subspace then the convex closure of ${\mathcal A}$
coincides with this parabolic subspace.

\begin{cor}\label{cor4-2}
Let $\Pi'$ be a polar space of rank $l$ and
let $m$ be a natural number such that the pair $l,m$ satisfies \eqref{eq4-1}.
If ${\mathcal S}$ is a convex subspace of  ${\mathfrak G}_{k}(\Pi)$ isomorphic to ${\mathfrak G}_{m}(\Pi')$
then there exists a $(k-m-1)$-dimensional singular subspace $N$ such that ${\mathcal S}=[N\rangle_{k}$.
\end{cor}

\begin{proof}
Let $f$ be a collineation of ${\mathfrak G}_{m}(\Pi')$ to ${\mathcal S}$
and let ${\mathcal A}$ be an apartment of ${\mathcal G}_{m}(\Pi')$.
By Theorem \ref{theorem4-2},
$f({\mathcal A})$ is an apartment in a parabolic subspace of ${\mathfrak G}_{k}(\Pi)$
or it is contained in a big star $[S\rangle_{k}$, $S\in {\mathcal G}_{k-1}(\Pi)$.
In the second case, the image of $\langle {\mathcal A}\rangle$
is also contained in $[S\rangle_{k}$. Since $[S\rangle_{k}$ is a polar space,
every point of $f(\langle {\mathcal A}\rangle)$ is collinear with at least one point on every line
contained in $f(\langle {\mathcal A}\rangle)$; but the same does not hold for $\langle {\mathcal A}\rangle$.
So,  there is a $(k-m-1)$-dimensional singular subspace $N$ such that
$f({\mathcal A})$ is an apartment in $[N\rangle_{k}$.
The subspace ${\mathcal S}$ contains the convex closure of $f({\mathcal A})$, i.e. $[N\rangle_{k}$.
Then $f^{-1}([N\rangle_{k})$ is a convex subspace of ${\mathfrak G}_{m}(\Pi')$ containing ${\mathcal A}$
which implies that ${\mathcal S}=[N\rangle_{k}$.
\end{proof}

\begin{rem}{\rm
We keep the notations of Corollary \ref{cor4-2}.
Every collineation of ${\mathfrak G}_{m}(\Pi')$ to $[N\rangle_{k}$ is induced by a collineation of $\Pi'$ to the polar space $[N\rangle_{k-m}$
except the case when our polar spaces are of type $\textsf{D}_4$ and $m=1$.
In this  case, every collineation of ${\mathfrak G}_{m}(\Pi')$ to $[N\rangle_{k}$ is induced by
a collineation of $\Pi'$ to $[N\rangle_{k-m}$ or a collineation of $\Pi'$ to one of the half-spin Grassmann spaces of $[N\rangle_{k-m}$.
See \cite[Section 4.6]{Pankov-book} for the details.
}\end{rem}

If ${\mathcal X}$ is a convex subspace  of ${\mathfrak G}_{k}(\Pi)$ isomorphic
to a polar space of rank $n-k\ge 3$
then one of the following possibilities is realized:
\begin{enumerate}
\item[$\bullet$] ${\mathcal X}$ is a big star of ${\mathcal G}_{k}(\Pi)$,
\item[$\bullet$] $n-k=3$ and ${\mathcal X}=[N,M]_{k}$ for some incident
$N\in {\mathcal G}_{k-2}(\Pi)$, $M\in {\mathcal G}_{k+2}(\Pi)$.
\end{enumerate}
This follows from Proposition \ref{prop3-1} and the fact that every polar space is the convex closure of any frame.

\subsection{Generalizations of Theorems \ref{theorem4-2} and  \ref{theorem4-3}}
By Proposition \ref{prop3-2}, embeddings of $PJ(l,m)$ in $\Gamma_{k}(\Pi)$ exist only in the case when $l-m\le n-k$.
If a such embedding sends tops of $PJ(l,m)$ to independent subsets of ${\mathfrak G}_{k}(\Pi)$ contained in tops then
$m\le k$. Indeed, every top of $PJ(l,m)$ consists of $m+2$ elements and eve\-ry top of ${\mathcal G}_{k}(\Pi)$
is a $(k+1)$-dimensional projective space. So, we get the following weak version of the condition \eqref{eq4-1}:
\begin{equation}\label{eq4-2}
0<m\le k,\;l\le n\;\mbox{ and }\;3\le l-m\le n-k.
\end{equation}

\begin{theorem}\label{theorem4-4}
Let $l,m$ be natural numbers satisfying \eqref{eq4-2}.
Let also $f$ be an embedding of $PJ(l,m)$ in $\Gamma_{k}(\Pi)$ transferring tops
to independent subsets of ${\mathfrak G}_{k}(\Pi)$ contained in tops.
Then there exists a $(k-m-1)$-dimensional singular subspace $N$ such that
the image of $f$ is contained in $[N\rangle_{k}$;
there is also an $l$-frame ${\mathcal B}\subset[N\rangle_{k-m}$
such that every element in the image of $f$
is spanned by elements of ${\mathcal B}$.
\end{theorem}

If $l-m=n-k$ then ${\mathcal B}$ is a frame of the polar space $[N\rangle_{k-m}$.
So, Theorem \ref{theorem4-3} is a partial case of Theorem \ref{theorem4-4}.

In the next section,
we prove Theorem \ref{theorem4-4} and give an independent proof of Theorem \ref{theorem4-3} based on a general result
concerning apartments in building Grassmannians \cite[Theorem 3.3]{Kasikova1}.
In the case when $l-m<n-k$, the images of all embeddings of $PJ(l,m)$ in $\Gamma_{k}(\Pi)$
are not apartments of building Grassmannians. By this reason,
\cite[Theorem 3.3]{Kasikova1} can not be used to prove Theorem \ref{theorem4-4}.

There is the following analog of Theorem \ref{theorem4-2}.

\begin{theorem}\label{theorem4-5}
Let $l,m$ be natural numbers satisfying \eqref{eq4-2}.
Let also ${\mathcal X}$ be a subset of ${\mathcal G}_{k}(\Pi)$ such that
$\Gamma({\mathcal X})$ is isomorphic to $PJ(l,m)$ and
every maximal clique of $\Gamma({\mathcal X})$ is an independent subset of ${\mathfrak G}_{k}(\Pi)$.
Suppose that one of the following conditions holds:
\begin{enumerate}
\item[{\rm (1)}]
$m+2>n-k$,
\item[{\rm (2)}]
$l-m\ge 4$ and ${\mathcal X}$ is not contained in any big star $[S\rangle_{k}$, $S\in {\mathcal G}_{k-1}(\Pi)$.
\end{enumerate}
Then there exists a $(k-m-1)$-dimensional singular subspace $N$ such that ${\mathcal X}$ is contained in
$[N\rangle_{k}$ and there is an $l$-frame ${\mathcal B}\subset[N\rangle_{k-m}$
such that every element of ${\mathcal X}$ is spanned by elements of ${\mathcal B}$.
\end{theorem}

\begin{rem}\label{rem4-2}{\rm
Suppose that $n-k=l-m$. Then (1) coincides with the condition (1) from Theorem \ref{theorem4-2}
and if $l-m=3$ then (1) holds for all $m\ge2$.
This means that Theorem \ref{theorem4-2} is a consequence of Theorem \ref{theorem4-5}
in the following cases:
\begin{enumerate}
\item[$\bullet$] $l-m\ge 4$,
\item[$\bullet$] $l-m=3$ and $m\ge 2$.
 \end{enumerate}
It was established above that Theorem \ref{theorem4-2} follows from Theorem \ref{theorem4-1} if $l-m=2$.
Therefore, it is sufficient to prove Theorem \ref{theorem4-2} only in the case when $l=4$, $m=1$
if Theorem \ref{theorem4-5} is proved.
}\end{rem}

By Propositions \ref{prop3-3}, every embedding of $PJ(l,1)$ in $\Gamma_{k}(\Pi)$ sends all maximal cliques of $PJ(l,1)$
to independent subsets of ${\mathfrak G}_{k}(\Pi)$ if $l\ge 4$. This implies the following.

\begin{cor}\label{cor4-3}
Suppose that $n-k+1\ge l$ and ${\mathcal X}$ is a subset of ${\mathcal G}_{k}(\Pi)$ such that
$\Gamma({\mathcal X})$ is isomorphic to $PJ(l,1)$.
If ${\mathcal X}$ is not contained in any big star $[S\rangle_{k}$, $S\in {\mathcal G}_{k-1}(\Pi)$
then there exists a $(k-2)$-dimensional singular subspace $N$ such that ${\mathcal X}$ is contained in
$[N\rangle_{k}$ and there is an $l$-frame ${\mathcal B}\subset[N\rangle_{k-1}$
such that every element of ${\mathcal X}$ is spanned by elements of ${\mathcal B}$.
\end{cor}

\section{Proof of Theorems \ref{theorem4-3} and \ref{theorem4-4}}
\subsection{}
Let $B=\{p_{1},\dots,p_{2l}\}$ be a frame of a certain rank $l$ polar space $\Pi'$.
For each $i\in \{0,\dots,l-1\}$ the associated apartment of
${\mathcal G}_i(\Pi')$ will be denoted by ${\mathcal A}_i$.

We say that a mapping $f:{\mathcal A}_{i}\to {\mathcal G}_{j}(\Pi)$
is {\it adjacency preserving} if
it sends any pair of adjacent elements of ${\mathcal A}_{i}$ to adjacent elements of ${\mathcal G}_{j}(\Pi)$.
An adjacency preserving mapping is not necessarily injective, but its restriction to every clique of $\Gamma({\mathcal A}_{i})$ is injective.
This mapping transfers maximal cliques of $\Gamma({\mathcal A}_{i})$ to subsets in maximal cliques of $\Gamma_{j}(\Pi)$,
but the images of distinct maximal cliques are not necessarily contained in distinct maximal cliques.

\begin{lemma}\label{lemma5-1}
Suppose that $0<i\le l-3$, $0<j\le n-3$ and  $f:{\mathcal A}_{i}\to {\mathcal G}_{j}(\Pi)$
is an adjacency preserving mapping which
sends stars to subsets of stars. Then big stars go to subsets of big stars;
moreover, the image of every big star is contained in a unique big star.
\end{lemma}

\begin{proof}
Let $S\in {\mathcal A}_{i-1}$.
The big star ${\mathcal B}_i(S)$ is the sum of all stars ${\mathcal S}_i(S,M)$, where $M\in {\mathcal A}_{l-1}$ contains $S$.
For every such $M$ there exist
$$g_{S}(M)\in {\mathcal G}_{j-1}(\Pi)\;\mbox{ and }\;h_{S}(M)\in {\mathcal G}_{n-1}(\Pi)$$
such that
$$f({\mathcal S}_i(S,M))\subset [g_{S}(M),h_{S}(M)]_j.$$
If $M$ and $M'$ are adjacent elements of ${\mathcal A}_{l-1}\cap[S\rangle_{l-1}$
then the intersection of the stars ${\mathcal S}_i(S,M)$ and ${\mathcal S}_i(S,M')$ contains at least two elements.
The same holds for the intersection of their images
(since the restriction of $f$ to every clique of $\Gamma({\mathcal A}_{i})$ is injective), i.e.
$$|[g_{S}(M),h_{S}(M)]_j\cap[g_{S}(M'),h_{S}(M')]_j|\ge 2.$$
The latter implies that
$$g_{S}(M)=g_{S}(M').$$
By connectedness (any $M,M'\in {\mathcal A}_{l-1}\cap[S\rangle_{l-1}$ can be connected by a path of $\Gamma_{l-1}(\Pi)$
contained in ${\mathcal A}_{l-1}\cap[S\rangle_{l-1}$),
this equality holds for all $M,M'\in {\mathcal A}_{l-1}$ containing $S$.
Therefore, the image of ${\mathcal B}_i(S)$ is in a big star of ${\mathcal G}_{j}(\Pi)$.
Since the intersection of two distinct big stars of ${\mathcal G}_{j}(\Pi)$ contains at most one element,
there is only one big star of ${\mathcal G}_{j}(\Pi)$ containing $f({\mathcal B}_i(S))$.
\end{proof}

\subsection{Proof of Theorem \ref{theorem4-4}}
Let $m$ be a natural number such that the pair $l,m$ satisfies \eqref{eq4-2}.
Let also $i\in \{0,\dots,m-1\}$ and let $f:{\mathcal A}_{m-i}\to {\mathcal G}_{k-i}(\Pi)$ be a mapping
satisfying the following condition:
\begin{enumerate}
\item[(T1)] the restriction of the mapping to every top is injective and
the image of every top is an independent subset contained in a top,
\end{enumerate}
Any two adjacent elements of ${\mathcal A}_{m-i}$ are contained in a top and (T1) implies that $f$ is adjacency preserving.
Hence its restriction to every clique of $\Gamma({\mathcal A}_{m-i})$ is injective.
Since the intersection of two distinct tops of ${\mathcal G}_{k}(\Pi)$ contains at most one element,
the image of every top  is contained in a unique top.
We say that the mapping $f$ is {\it special} if it satisfies the following additional condition:
\begin{enumerate}
\item[(T2)] the images of adjacent tops\footnote{We say that two tops of ${\mathcal A}_{j}$,
$j\in \{1,\dots,l-2\}$ are {\it adjacent} if the associated elements of ${\mathcal A}_{j+1}$ are adjacent.}
are contained in distinct tops.
\end{enumerate}

\begin{lemma}\label{lemma5-2}
Every special mapping $f:{\mathcal A}_{m-i}\to {\mathcal G}_{k-i}(\Pi)$ transfers stars to subsets of stars.
\end{lemma}

\begin{proof}
For every star ${\mathcal S}\subset {\mathcal A}_{m-i}$ there is a top ${\mathcal T}\subset {\mathcal A}_{m-i}$
intersecting ${\mathcal S}$ precisely in two elements.
By $({\rm T}1)$, $f({\mathcal T})$ is in a top.
Since the intersection of two distinct tops of ${\mathcal G}_{k-i}(\Pi)$ contains at most one element
and $|f({\mathcal S})\cap f({\mathcal T})|=2$, $f({\mathcal S})$ is in a star
or in the top containing $f({\mathcal T})$.

Suppose that $f({\mathcal S})$ is not in a star.
Then it is in the top containing $f({\mathcal T})$. We take any top ${\mathcal T}'\subset {\mathcal A}_{m-i}$ adjacent with ${\mathcal T}$ and such that $|{\mathcal S}\cap {\mathcal T}'|=2$.
As above, $f({\mathcal S})$ is in a star or in the top containing  $f({\mathcal T}')$.
Since $f({\mathcal S})$ is not in a star,  it is in the top containing $f({\mathcal T}')$.
By $({\rm T}2)$, the tops containing $f({\mathcal T})$ and $f({\mathcal T}')$ are distinct and their intersection contains at most one element.
So, $f({\mathcal S})$ can not be contained in a top.
\end{proof}

Let $f:{\mathcal A}_{m}\to {\mathcal G}_{k}(\Pi)$ be a special mapping.
By Lemmas \ref{lemma5-1} and \ref{lemma5-2}, it transfers big stars to subsets of big stars and the image of every big star
is contained in a unique big star.
Therefore, there is a mapping
$f_{m-1}:{\mathcal A}_{m-1}\to {\mathcal G}_{k-1}(\Pi)$
such that
$$f({\mathcal B}_{m}(S))\subset [f_{m-1}(S)\rangle_{k}\;\;\;\;\;\forall\;S\in {\mathcal A}_{m-1}.$$
If $m>1$ then
$$f_{m-1}({\mathcal T}_{m-1}(U))\subset \langle f(U)]_{k-1}\;\;\;\;\;\forall\;U\in {\mathcal A}_m,$$
i.e. $f_{m-1}$ transfers tops to subsets of tops.
If $m=1$ then $f_{0}$ sends any pair of distinct collinear points of $B$ to elements of ${\mathcal G}_{k-1}(\Pi)$
which are adjacent or coincident.

\begin{lemma}\label{lemma5-3}
If  $m>1$ then $f_{m-1}$ is special.
If $k=1$ then $f_{0}$ transfers any pair of distinct collinear points of $B$ to
adjacent elements of ${\mathcal G}_{k-1}(\Pi)$.
\end{lemma}

\begin{proof}
Let $m>1$ and $U\in {\mathcal A}_m$. Denote by $S_1,\dots,S_{m+1}$
the elements of the top ${\mathcal T}_{m-1}(U)$.
We take any point $p_{j}\in B\setminus U$.
Then $\langle p_{j}, U \rangle$ belongs to ${\mathcal A}_{m+1}$ and
the top ${\mathcal T}_{m}(\langle p_{j}, U \rangle)$ consists of
$$S'_{i}:=\langle p_{j}, S_i \rangle,\;\;\;i=1,\dots,m+1$$
and $U$.
By (T1),
$$f(U),\,f(S'_{1}),\,\dots,\,f(S'_{m+1})$$
are distinct and form an independent subset in a certain top of ${\mathcal G}_{k}(\Pi)$.
Then
$$f(U)\cap f(S'_{1}),\dots, f(U)\cap f(S'_{m+1})$$
are distinct and  form an independent subset in the top $\langle f(U)]_{k-1}$.
It is clear that
$$f_{k-1}(S_i)=f(U)\cap f(S'_{i}).$$
Therefore, $f_{m-1}$ satisfies (T1) and
the same arguments show that $f_{0}$ maps any pair of distinct collinear points of $B$ to
adjacent elements of ${\mathcal G}_{k-1}(\Pi)$ if $m=1$.
Since $f_{m-1}$ is induced by $f$ and $f$ is adjacency preserving,
$f_{m-1}$ satisfies (T2).
\end{proof}

\begin{rem}{\rm
We are not be able to show that $f_{m-1}$ is an embedding of $PJ(l,m-1)$ in $\Gamma_{k-1}(\Pi)$ if
$f$ is an embedding of $PJ(l,m)$ in $\Gamma_{k}(\Pi)$.
By this reason, we consider the class of special mappings instead of embeddings.
}\end{rem}

Step by step, we construct a sequence of mappings
$$f_{i}:{\mathcal A}_{i}\to {\mathcal G}_{k-m+i}(\Pi),\;\;\;i=m,m-1,\dots,0$$
such that $f_m=f$ and $f_{i}$ is special if $i\ge 1$. Moreover, we have
$$f_i({\mathcal B}_{i}(S))\subset [f_{i-1}(S)\rangle_{k-m+i}\;\;\;\;\;\forall\;S\in {\mathcal A}_{i-1},$$
$$f_{i-1}({\mathcal T}_{i-1}(U))\subset \langle f_i(U)]_{k-m+i-1}\;\;\;\;\;\forall\;U\in {\mathcal A}_i$$
for every $i\ge 1$.
The latter implies that if $X\in {\mathcal A}_{i}$ and $Y\in {\mathcal A}_{j}$, $i,j\le m$ are incident
then $f_{i}(X)$ and $f_{j}(Y)$ are incident.

For every $i\in \{1,\dots,2l\}$ we define $Q_{i}:=f_{0}(p_{i})$.
Lemma \ref{lemma5-3} states that $Q_{i}$ and $Q_{j}$, $j\ne i,\sigma(i)$ are adjacent elements of ${\mathcal G}_{k-m}(\Pi)$.

\begin{lemma}\label{lemma5-4}
If $p_{i_{1}},\dots,p_{i_{j+1}}$ span an element of ${\mathcal A}_{j}$, $j\le m$ then
$$f_{j}(\langle p_{i_{1}},\dots,p_{i_{j+1}} \rangle)=\langle Q_{i_{1}},\dots,Q_{i_{j+1}} \rangle.$$
\end{lemma}

\begin{proof}
The case $j=0$ is trivial.
Suppose that $j>0$ and prove the statement by induction.
By inductive hypothesis, $f_{j-1}$ transfers
$$\langle p_{i_{1}},\dots,p_{i_{j-1}},p_{i_{j}}\rangle\;\mbox{ and }\;\langle p_{i_{1}},\dots,p_{i_{j-1}},p_{i_{j+1}}\rangle$$
to
$$\langle Q_{i_{1}},\dots,Q_{i_{j-1}},Q_{i_{j}}\rangle\;\mbox{ and }\;\langle Q_{i_{1}},\dots,Q_{i_{j-1}},Q_{i_{j+1}}\rangle,$$
respectively. These are adjacent elements of ${\mathcal G}_{k-m+j-1}(\Pi)$
contained in $$f_{j}(\langle p_{i_{1}},\dots,p_{i_{j+1}} \rangle)\in{\mathcal G}_{k-m+j}(\Pi) $$
and we get the required equality.
\end{proof}

Let $I$  be a subset of $\{1,\dots,2l\}$ consisting of $l$ elements and satisfying $I\cap\sigma(I)=\emptyset$.
Then $\{Q_i\}_{i\in I}$ is a clique in $\Gamma_{k-m}(\Pi)$.
If this clique is contained in a top then for any distinct $u,v,w\in I$ we have
$$f_{1}(p_{u}p_{v})=\langle Q_{u},Q_{v} \rangle=\langle Q_{u},Q_{w} \rangle=f_{1}(p_{u}p_{w})$$
which is impossible.
Therefore, $$N:=\bigcap_{i\in I}Q_i\in {\mathcal G}_{k-m-1}(\Pi).$$
For every $j\in \{1,\dots,2l\}\setminus I$ there exist $i_{1},i_{2}\in I$
such that $Q_{j},Q_{i_{1}},Q_{i_{2}}$ form a clique in $\Gamma_{k-m}(\Pi)$.
As above, we show that this clique can not be contained in a top. Then $N=Q_{i_{1}}\cap Q_{i_{2}}$
is contained in $Q_{j}$.

So, every $Q_{i}$ belongs to $[N\rangle_{k-m}$ and
$Q_i,Q_j$ are collinear points of the polar space $[N\rangle_{k-m}$ if $j\ne i,\sigma(i)$.

Suppose that $f$ is an embedding of $PJ(l,m)$ in  $\Gamma_{k}(\Pi)$.
For every $i\in \{1,\dots,2l\}$ we can choose $i_{1},\dots,i_{m}\in \{1,\dots,2l\}\setminus\{i,\sigma(i)\}$
such that $$\{i_{1},\dots,i_{m}\}\cap \{\sigma(i_{1}),\dots,\sigma(i_{m})\}=\emptyset.$$
Then
$$\langle p_{i_{1}},\dots,p_{i_{m}},p_{i}\rangle\;\mbox{ and }\;\langle p_{i_{1}},\dots,p_{i_{m}},p_{\sigma(i)}\rangle$$
are non-adjacent elements of ${\mathcal A}_{m}$ and
their images
$$\langle Q_{i_{1}},\dots,Q_{i_{m}},Q_{i}\rangle\;\mbox{ and }\;\langle Q_{i_{1}},\dots,Q_{i_{m}},Q_{\sigma(i)}\rangle$$
are non-adjacent elements of ${\mathcal G}_{k}(\Pi)$.
Hence $Q_{i}$ are $Q_{\sigma(i)}$ are non-collinear points of the polar space $[N\rangle_{k-m}$.
The rank of this polar space is $n-k+m\ge l$
and $Q_{1},\dots,Q_{2l}$ form  an $l$-frame in $[N\rangle_{k-m}$.
By Lemma \ref{lemma5-4}, every element in the image of $f$ is spanned by elements of this $l$-frame.

\subsection{Proof of Theorem \ref{theorem4-3}}
Suppose that $f:{\mathcal A}_{m}\to {\mathcal G}_{k}(\Pi)$ is an embedding of $PJ(l,m)$ in $\Gamma_{k}(\Pi)$ and
the condition \eqref{eq4-1} holds.
The image of $f$ will be denoted by ${\mathcal X}$.
We assume that $f$ transfers tops to independent subsets of ${\mathfrak G}_{k}(\Pi)$
contained in tops.  Subsets of type
$${\mathcal X}\cap f({\mathcal T}_{m}(U)),\;\;U\in {\mathcal A}_{m+1}$$
are said to be {\it tops} of ${\mathcal X}$.
By Lemma \ref{lemma5-2}, $f$ sends stars to subsets of stars
and Lemma \ref{lemma5-1} implies that big stars go to subsets of big stars\footnote{In the case when $l-m\ge 4$,
this follows from Proposition \ref{prop3-2}.}.

By \cite[Theorem 3.3]{Kasikova1}, ${\mathcal X}$ is an apartment in a parabolic subspace of ${\mathfrak G}_{k}(\Pi)$ if
the following holds:
\begin{enumerate}
\item[(A)] for every $S\in {\mathcal X}$ there exist a $(k-m-1)$-dimensional singular subspace $N_{S}$
and an apartment ${\mathcal A}_{S}$ in the parabolic subspace $[N_{S}\rangle_{k}$
such that $S\in {\mathcal A}_{S}$ and ${\mathcal X}(S)={\mathcal A}_{S}(S)$\footnote{For any subset ${\mathcal X}\subset {\mathcal G}_{k}(\Pi)$
and any $S\in {\mathcal X}$ we denote by ${\mathcal X}(S)$ the set consisting of $S$ and all elements of
${\mathcal X}$ adjacent with $S$.}.
\end{enumerate}
There is also a condition concerning three elements of ${\mathcal X}$ contained in a plane,
but in our case it is obvious.

We take any top
$${\mathcal T}={\mathcal X}\cap \langle U_{\mathcal T}]_{k},\;\;\;U_{\mathcal T}\in {\mathcal G}_{k+1}(\Pi).$$
Let $S_{1},\dots,S_{m+2}$ be the elements of ${\mathcal T}$. Since they form an independent subset of $\langle U]_{k}$,
the singular subspace
$$N_{\mathcal T}:=S_{1}\cap\dots\cap S_{m+2}$$
is $(k-m-1)$-dimensional and
$$X_{i}:=\bigcap_{j\ne i} S_{j},\;\;\;i=1,\dots,m+2$$
form a base in the projective space $[N_{\mathcal T},U_{\mathcal T}]_{k-m}$.
Denote this base by ${\mathcal B}_{\mathcal T}$.

\begin{lemma}\label{lemma5-5}
Let ${\mathcal T}$ and ${\mathcal T}'$ be distinct tops of ${\mathcal X}$ containing $S$.
Then for every $Y\in {\mathcal T}\setminus \{S\}$ there is unique $Y'\in {\mathcal T}'\setminus \{S\}$
such that
$$Y\cap S=Y'\cap S.$$
\end{lemma}

\begin{proof}
Consider $S'=f^{-1}(S)\in {\mathcal A}_m$ and $i,j\in \{1,\dots,2l\}$ such that
$$f^{-1}({\mathcal T})={\mathcal T}_{m}(\langle S', p_{i}\rangle)\;\mbox{ and }\;
f^{-1}({\mathcal T}')={\mathcal T}_{m}(\langle S', p_{j}\rangle).$$
If $Y\in {\mathcal T}\setminus \{S\}$ then $f^{-1}(Y)$ is spanned by $p_{i}$ and some $K\in {\mathcal A}_{m-1}$
contained in $S'$.
Then
$$Y':=f(\langle K, p_{j}\rangle)\in {\mathcal T}'$$
is as required.
Indeed, $f^{-1}(S)$, $f^{-1}(Y)$, $f^{-1}(Y')$
belong to the big star ${\mathcal B}_{m}(K)$;
since $f$ maps big stars to subsets of big stars, $S$, $Y$, $Y'$ are contained in a big star which implies the required equality.
\end{proof}

Let $S\in {\mathcal X}$.
Then ${\mathcal X}(S)$ coincides with the sum of all tops ${\mathcal T}\subset {\mathcal X}$ containing $S$.
If a top ${\mathcal T}\subset {\mathcal X}$ contains $S$ then ${\mathcal B}_{\mathcal T}$ consists of
$m+1$ elements contained in $S$ and one element not in $S$, this element will be denoted by $X_{\mathcal T}$.
Lemma \ref{lemma5-5} implies that
for any distinct tops ${\mathcal T},{\mathcal T}'\subset {\mathcal X}$ containing $S$
the following assertions are fulfilled:
\begin{enumerate}
\item[(1)] $N_{\mathcal T}=N_{{\mathcal T}'}$,
i.e. $N_{\mathcal T}$ does not depend on ${\mathcal T}$ and we denote this subspace by $N_{S}$;
\item[(2)]
${\mathcal B}_{\mathcal T}\cap {\mathcal B}_{{\mathcal T}'}$ consists of $m+1$ elements contained in $S$
and ${\mathcal B}_{\mathcal T}\setminus{\mathcal B}_{{\mathcal T}'}=\{X_{\mathcal T}\}$.
\end{enumerate}
Let ${\mathcal B}_S$ be the sum of all ${\mathcal B}_{\mathcal T}$ such that $S\in {\mathcal T}$.
Since there are precisely  $2(l-m-1)$ distinct tops of ${\mathcal X}$ containing $S$, the statement (2) shows that
${\mathcal B}_S$ is formed by $m+1$ elements contained in $S$ and $2(l-m-1)$ elements not in $S$.

From this moment, we will consider ${\mathcal B}_S$ as a subset in the polar space $[N_{S}\rangle_{k-m}$.
Every $X\in {\mathcal B}_S$ satisfying $X\subset S$ is collinear with all other points of ${\mathcal B}_S$.

\begin{lemma}\label{lemma5-6}
Let ${\mathcal T}$ and ${\mathcal T}'$ be distinct tops of ${\mathcal X}$ containing $S$.
Then $X_{\mathcal T}$ and $X_{{\mathcal T}'}$
are collinear if and only if the tops $f^{-1}({\mathcal T})$ and $f^{-1}({\mathcal T}')$ are adjacent.
\end{lemma}

\begin{proof}
Let $Y\in {\mathcal T}\setminus \{S\}$. By Lemma \ref{lemma5-5},
there is unique $Y'\in {\mathcal T}'\setminus \{S\}$ such that
$$Y\cap S=Y'\cap S.$$
The tops
$f^{-1}({\mathcal T})$ and $f^{-1}({\mathcal T}')$ are adjacent if and only if $f^{-1}(Y)$ and $f^{-1}(Y')$ are adjacent,
i.e. $Y$ and $Y'$ are adjacent (see the proof of Lemma \ref{lemma5-5}).
The latter is equivalent to the fact that
the elements of ${\mathcal G}_{k+1}(\Pi)$ associated with ${\mathcal T}$ and ${\mathcal T}'$ are adjacent
which gives the claim.
\end{proof}

\begin{lemma}\label{lemma5-7}
For every $X\in {\mathcal B}_S$ satisfying $X\not\subset S$
there is unique point of ${\mathcal B}_S$ non-collinear with $X$.
\end{lemma}

\begin{proof}
There is unique top ${\mathcal T}\subset {\mathcal X}$ containing $S$
and such that $X=X_{\mathcal T}$.
Suppose that $f^{-1}({\mathcal T})={\mathcal T}_{m}(U)$.
Then $U$ is spanned by $f^{-1}(S)$ and a point $p_{i}$.
Let $U'$ be the subspace spanned by $f^{-1}(S)$ and  $p_{\sigma(i)}$.
Then
$${\mathcal T}':=f({\mathcal T}_{m}(U'))$$
is unique top of ${\mathcal X}$ containing $S$ and whose preimage is not adjacent with $f^{-1}({\mathcal T})$.
Lemma \ref{lemma5-6} guarantees that $X_{{\mathcal T}'}$ is unique
point of ${\mathcal B}_S$ non-collinear with $X$.
\end{proof}

By Lemma \ref{lemma2-2}, ${\mathcal B}_S$ can be extended to a frame of the polar space $[N_{S}\rangle_{k-m}$.
If ${\mathcal A}_{S}$ is the associated apartment of $[N_{S}\rangle_{k}$ then
$S\in{\mathcal A}_{S}$ and  ${\mathcal X}(S)={\mathcal A}_{S}(S)$.

\section{Proof of Theorems \ref{theorem4-2} and \ref{theorem4-5}}
\setcounter{equation}{0}
\subsection{Proof of Theorem \ref{theorem4-5}}
Let $l$ and $m$ be natural numbers satisfying \eqref{eq4-2}.

\begin{lemma}\label{lemma6-1}
Let $f$ be an embedding of $PJ(l,m)$ in $\Gamma_{k}(\Pi)$ transferring stars
to subsets of star.
If the image of a certain top is contained in a star
then the image of $f$ is in a big star.
\end{lemma}

By Proposition \ref{prop3-1}, we have the following.

\begin{cor}\label{cor6-1}
If $l-m\ge 4$ and an embedding of $PJ(l,m)$ in $\Gamma_{k}(\Pi)$
sends a certain top to a subset of a star
then the image of this embedding is contained in a big star.
\end{cor}

\begin{proof}[Proof of Lemma \ref{lemma6-1}]
Let $B$ be a frame of a rank $l$ polar space $\Pi'$.
For every $i\in\{0,\dots,l-1\}$ we denote by ${\mathcal A}_{i}$ the associated apartment of ${\mathcal G}_{i}(\Pi')$.
Let also $f:{\mathcal A}_{m}\to {\mathcal G}_{k}(\Pi)$ be an embedding of $PJ(l,m)$ in $\Gamma_{k}(\Pi)$
transferring stars to subsets of stars. By Lemma \ref{lemma5-1}, big stars go to subsets of big stars.
Suppose that the image of ${\mathcal T}_{m}(T)$, $T\in{\mathcal A}_{m+1}$ is contained in a star. Then
$$f({\mathcal T}_{m}(T))\subset [S\rangle_{k}\;\mbox{ for some }\; S\in {\mathcal G}_{k-1}(\Pi).$$
If $U\in {\mathcal A}_{m-1}$ is contained in $T$ then the intersection of the associated big star
${\mathcal B}_{m}(U)$ with ${\mathcal T}_{m}(T)$ consists of $2$ elements.
Since big stars go to subsets of big stars
and the intersection of two distinct big stars of ${\mathcal G}_{k}(\Pi)$ contains at most one element,
the inclusion
$$f({\mathcal B}_{m}(U))\subset [S\rangle_{k}$$
holds for every $U\in {\mathcal A}_{m-1}$ contained in $T$.

Let $T'$ be an element of ${\mathcal A}_{m+1}$ adjacent with $T$.
Then $X:=T\cap T'\in {\mathcal A}_{m}$ and
$${\mathcal T}_{m}(T)\cap {\mathcal T}_{m}(T')=\{X\}.$$
If $Y\in {\mathcal T}_{m}(T')\setminus\{X\}$ then $X\cap Y$ is an element of ${\mathcal A}_{m-1}$ contained in $T$
and we have
$$f(Y)\in f({\mathcal B}_{m}(X\cap Y))\subset [S\rangle_{k}.$$
Therefore,
$$f({\mathcal T}_{m}(T'))\subset [S\rangle_{k}.$$
By connectedness, this inclusion holds for every $T'\in {\mathcal A}_{m+1}$.
So, the image of $f$ is contained in the big star $[S\rangle_{k}$.
\end{proof}

Let ${\mathcal X}$ be a subset of ${\mathcal G}_{k}(\Pi)$ such that
$\Gamma({\mathcal X})$ is isomorphic to $PJ(l,m)$ and every maximal clique of $\Gamma({\mathcal X})$ is an independent subset of
${\mathfrak G}_{k}(\Pi)$. Then ${\mathcal X}$ is the image of an embedding $f$ of $PJ(l,m)$ in $\Gamma_{k}(\Pi)$
which transfers maximal cliques of $PJ(l,m)$ to independent subsets of ${\mathfrak G}_{k}(\Pi)$.
The maximal cliques of $PJ(l,m)$ are stars and tops;
they contain $l-m$ and $m+2$ elements, respectively.
Their images are independent subsets in maximal cliques of $\Gamma_{k}(\Pi)$ --- stars and tops
which are projective spaces of dimension $n-k-1$ and $k+1$, respectively.
By Theorem \ref{theorem4-4}, it sufficient to show that
$f$ sends tops to subsets of tops.
If $m+2>n-k$ (the first condition of Theorem \ref{theorem4-5})
then the images of tops can not be independent subsets of stars; thus the image of every top is contained in a top.
In the case when $l-m\ge 4$ and ${\mathcal X}$ is not contained in any big star of ${\mathcal G}_{k}(\Pi)$
(the second condition of Theorem \ref{theorem4-5}), Corollary \ref{cor6-1} gives the claim.

\subsection{Proof of Theorem \ref{theorem4-2}}
By Remark \ref{rem4-2}, we restrict ourself to the case when $l=4$ and $m=1$.

Let $B=\{p_{1},\dots,p_{8}\}$ be a frame of a certain rank $4$ polar space $\Pi'$.
For every $i\in \{0,\dots,3\}$ we denote by ${\mathcal A}_{i}$ the associated apartment of ${\mathcal G}_{i}(\Pi')$.
There is unique presentation of ${\mathcal A}_3$ as the sum of two disjoint subsets
${\mathcal A}_{+}$ and ${\mathcal A}_{-}$ satisfying the following conditions:
\begin{enumerate}
\item[$\bullet$] $\dim (X\cap Y)$ belongs to $\{3,1,-1\}$ if $X,Y\in {\mathcal A}_{\delta}$, $\delta\in \{+,-\}$,
\item[$\bullet$] $\dim (X\cap Y)$ belongs to $\{2,0\}$ if $X\in {\mathcal A}_{\delta}$ and $Y\in {\mathcal A}_{-\delta}$;
\end{enumerate}
here $-\delta$ is the complement of $\delta$ in $\{+,-\}$.
Denote by $\Gamma_{\delta}$, $\delta\in \{+,-\}$ the graph whose vertex set is
${\mathcal A}_{\delta}$ and whose edges are pairs $X,Y\in{\mathcal A}_{\delta}$ satisfying
$$\dim (X\cap Y)=1.$$
Note that the graphs $\Gamma(B)=PJ(4,0)$ and $\Gamma_{\delta}$, $\delta\in \{+,-\}$ both are isomorphic to the $4$-dimensional half-cube graph $\frac{1}{2}H_4$.

Let $h:B\to {\mathcal A}_{-\delta}$ be an isomorphism between $\Gamma(B)$ and $\Gamma_{-\delta}$.
Then
$$h(p_i)\cap h(p_j)\in {\mathcal A}_{1}\;\mbox{ if }\;j\ne i,\sigma(i).$$
Denote by $g$ the transformation of ${\mathcal A}_{1}$
sending every line $p_{i}p_{j}$, $j\ne i,\sigma(i)$ to the intersection of $h(p_{i})$ and $h(p_{j})$.
This is an automorphism of $\Gamma({\mathcal A}_1)=PJ(4,1)$.
We will use the following fact:
{\it $g$ transfers every top to a star ${\mathcal S}_{1}(p_{u},U)$
such that $U\in {\mathcal A}_{\delta}$.}
\begin{proof}
Any top of ${\mathcal A}_{1}$ consists of three lines $p_{i}p_{j}$, $p_{i}p_{l}$ and $p_{j}p_{l}$.
Then $$h(p_{i}),h(p_{j}),h(p_{l})$$ are mutually adjacent elements of ${\mathcal A}_{-\delta}$.
By  \cite[Lemma 4.11]{Pankov-book},
$$h(p_{i})\cap h(p_{j})\cap h(p_{l})=g(p_{i}p_{j})\cap g(p_{i}p_{l})\cap g(p_{j}p_{l})$$
is a point and the lines
$$h(p_{i})\cap h(p_{j})=g(p_{i}p_{j}),\;h(p_{i})\cap h(p_{l})=g(p_{i}p_{l}),\;h(p_{j})\cap h(p_{l})=g(p_{j}p_{l})$$
span an element of ${\mathcal A}_{\delta}$.
\end{proof}

\begin{lemma}\label{lemma6-2}
Let $f$ be an embedding of $PJ(4,1)$ in $\Gamma_{k}(\Pi)$ such that the image of a certain star is contained in a top.
Then there exists an automorphism $g$ of $PJ(4,1)$
such that $fg$ maps tops to subsets of tops.
\end{lemma}

\begin{proof}
Suppose that $f:{\mathcal A}_{1}\to {\mathcal G}_{k}(\Pi)$ is an embedding of $PJ(4,1)$ in $\Gamma_{k}(\Pi)$
such that the image of a certain star ${\mathcal S}_{1}(p_{i}, U)$ is contained in a top $\langle N]_{k}$.
Our first step is to show that $f$ transfers every star ${\mathcal S}\subset{\mathcal A}_1$ formed by lines contained in $U$ to a subset of a top.

There are precisely $6$ elements of ${\mathcal A}_{1}$ contained in $U$ ---
$3$ lines belong to the star ${\mathcal S}_{1}(p_{i}, U)$ and the remaining $3$ lines form a top.
 Let $X_{1},X_{2},X_{3}$ and $Y_{1},Y_{2},Y_{3}$ be the elements of the star and the top, respectively.
Every element of the star is adjacent with precisely two elements of the top.
Assume that $X_{u}$ and $Y_{v}$ are adjacent if $u\ne v$.
Then $f(X_{u})$ is adjacent with $f(Y_{v})$ if and only if $u\ne v$.
By Proposition \ref{prop3-3}, $f(X_{1}),f(X_{2}),f(X_{3})$ form an independent subset in $\langle N]_{k}$. Hence
$$U_{1}:=f(X_{2})\cap f(X_{3}),\;U_{2}:=f(X_{1})\cap f(X_{3}),\;U_{3}:=f(X_{1})\cap f(X_{2})$$
are distinct elements of ${\mathcal G}_{k-1}(\Pi)$.
If $f(Y_{j})$ does not contain $U_j$ then it intersects $f(X_{u})$ and $f(X_{v})$, $u,v\ne j$ in
distinct elements of ${\mathcal G}_{k-1}(\Pi)$. This implies the inclusion
$$f(Y_{j})\subset \langle f(X_{u}),f(X_{v})\rangle=N.$$
Since $f(X_{j})\subset N$,
$f(Y_{j})$ is adjacent with $f(X_{j})$ which is impossible.
So, every $f(Y_{j})$ contains $U_j$.

Consider any star ${\mathcal S}\subset{\mathcal A}_1$ formed by lines contained in $U$
and distinct from the star ${\mathcal S}_{1}(p_{i}, U)$.
It consists of $X_{j}$ and $Y_{u},Y_{v}$ such that $u,v\ne j$. Then
$$f(X_{j})=\langle U_{u}, U_{v}\rangle\subset \langle f(Y_{u}),f(Y_{v})\rangle$$
which means that $f({\mathcal S})$ is contained in a top.

Let us take any $U'\in {\mathcal A}_{3}$ such that
$$\dim(U\cap U')=1$$
and any point $p_j$ on the line $L:=U\cap U'$.
The stars ${\mathcal S}_{1}(p_j, U)$ and ${\mathcal S}_{1}(p_j, U')$
both contain $L$.
Let $L_{1},L_{2}$ and $L'_{1},L'_{2}$ be the remaining elements of ${\mathcal S}_{1}(p_j, U)$ and ${\mathcal S}_{1}(p_j, U')$, respectively.
Every elements of ${\mathcal S}_{1}(p_j, U)\setminus\{L\}$ is adjacent with precisely one element of ${\mathcal S}_{1}(p_j, U')\setminus\{L\}$.
We suppose that $L_{i}$ is adjacent with $L'_{j}$ if $i=j$.
Then $f(L_{i})$ and $f(L'_{j})$ are adjacent if and only if $i=j$.
It was established above that $f({\mathcal S}_{1}(p_j, U))$ is contained in a top $\langle N']_k$ and
Proposition \ref{prop3-3} guarantees that $f(L),f(L_{1}),f(L_{2})$ form an independent subset in $\langle N']_k$. Hence
$$U_{1}:=f(L)\cap f(L_{1})\;\mbox{ and }\;U_{2}:=f(L)\cap f(L_{2})$$
are distinct elements of ${\mathcal G}_{k-1}(\Pi)$.

Since $f(L'_{1})$ is adjacent with $f(L),f(L_{1})$ and non-adjacent with $f(L_{2})$, we have
$U_{1}\subset f(L'_{1})$. Indeed,
if this fails then $f(L'_{1})$ intersects $f(L)$ and $f(L_{1})$ in distinct elements of ${\mathcal G}_{k-1}(\Pi)$
which implies the inclusion
$$f(L'_{1})\subset \langle f(L), f(L_{1})\rangle=N';$$
since $f(L_{2})\subset N'$, $f(L'_{1})$ is adjacent with $f(L_{2})$
which is impossible.
By the same reason, $U_{2}\subset f(L'_{2})$.
Therefore,
$$f(L)=\langle U_1,U_2 \rangle\subset \langle f(L'_{1}), f(L'_{2})\rangle $$
which means that the image of ${\mathcal S}_{1}(p_j, U')$ is contained in a top.

As above, we establish that the stars formed by elements  of ${\mathcal A}_1$ contained in $U'$
go to subsets of tops.

For one of $\delta\in\{+,-\}$ we have $U\in {\mathcal A}_{\delta}$.
Since the graph $\Gamma_{\delta}$ is connected,
the image of every star
\begin{equation}\label{eq6-1}
{\mathcal S}_{1}(p_u,U'),\;\;U'\in {\mathcal A}_{\delta}
\end{equation}
is contained in a top.
Let $g$ be an automorphism of $\Gamma({\mathcal A}_{1})$ transferring every top to a star of type \eqref{eq6-1}
(such automorphism was constructed above).
Then $fg$ is an embedding of $PJ(4,1)$ in $\Gamma_{k}(\Pi)$
sending tops to subset of tops.
\end{proof}

Let ${\mathcal X}$ be a subset of ${\mathcal G}_{k}(\Pi)$ such that
$\Gamma({\mathcal X})$ is isomorphic to $PJ(4,1)$.
We assume that ${\mathcal X}$ is not contained in any big star of ${\mathcal G}_{k}(\Pi)$.
Then ${\mathcal X}$ is the image of an embedding $f$ of $PJ(4,1)$ in $\Gamma_{k}(\Pi)$.
By Proposition \ref{prop3-3},
$f$ transfers maximal cliques of $PJ(l,m)$ to independent subsets of ${\mathfrak G}_{k}(\Pi)$.
If $f$ sends stars to subsets of stars then, by Lemma \ref{lemma6-1},
tops go to subsets of tops (since the image of $f$ is not contained in any big star)
and we apply Theorem \ref{theorem4-3}.
If the image of a certain star is contained in a top then Lemma  \ref{lemma6-2} implies the existence of
an automorphism $g$ of $J(4,1)$ such that $fg$ maps tops to subsets of tops.
We apply Theorem \ref{theorem4-3} to $fg$ and get the claim, since
the image of $fg$ is ${\mathcal X}$.


\begin{thebibliography}{999}
\bibitem{BB}
R. J. Blok, A. E. Brouwer,
Spanning point-line geometries in buildings of spherical type, J. Geom. 62(1998), 26--35.

\bibitem{CS}
B.N. Cooperstein, E.E. Shult,
Frames and bases of Lie incidence geometries, J. Geom.
60(1997), 17--46.

\bibitem{CKS}
B.N. Cooperstein, A. Kasikova, E.E. Shult,
Witt-type Theorems for Grassmannians and Lie Incidence Geometries,
Adv. Geom. 5(2005), 15--36.

\bibitem{Kasikova1} A. Kasikova,
Characterization of some subgraphs of point-collinearity graphs of building geometries,
European J. Combin. 28(2007), 1493--1529.

\bibitem{Kasikova2} A. Kasikova,
Convex subspace closure of the point shadow of an
apartment of a spherical building,
Adv. Geom. 9(2009), 541-561.

\bibitem{Kasikova3} A. Kasikova,
Characterization of some subgraphs of point-collinearity graphs of building geometries II, Adv. Geom. 9 (2009) 45--84.


\bibitem{Pankov-book}
Pankov M., Grassmannians of classical buildings, Algebra and
Discrete Math. Series 2, World Scientific, 2010.

\bibitem{Pankov1}
M. Pankov, Isometric embeddings of Johnson graphs in Grassmann graphs,
J. Algebraic Combin. 33(2011), 555--570.

\bibitem{Pankov2}
M. Pankov, Metric characterization of apartments in dual polar spaces,
J. Combin. Theory Ser. A 118(2011), 1313--1321.

\bibitem{Pasini}
Pasini A., Diagram geometries, Clarendon press, Oxford 1994.


\bibitem{Tits} J. Tits,
Buildings of spherical type and finite BN-pairs,
Lecture Notes in Mathematics 386,
Springer-Verlag, Berlin-New York, 1974.


\end{thebibliography}
\end{document}